\documentclass[a4paper]{amsart}
\usepackage{amsmath,amsfonts,amssymb,amsthm}
\usepackage{enumerate}
\usepackage{color}
\newtheorem{theorem}{Theorem}[section]
\newtheorem{lemma}[theorem]{Lemma}

\newtheorem{corollary}[theorem]{Corollary}
\theoremstyle{remark}

\theoremstyle{definition}
\newtheorem{definition}[theorem]{Definition}

\newcommand{\R}{\mathbb{R}}
\newcommand{\Z}{\mathbb{Z}}
\newcommand{\C}{\mathbb{C}}
\newcommand{\N}{\mathbb{N}}

\newcommand{\re}{\mathrm{Re}}
\newcommand{\im}{\mathrm{Im}}

\title{Applying the Resonance Method to $\re\left(e^{-i\theta}\log\zeta(\sigma+it)\right)$}
\author{Mikko Jaskari}
\address{Department of Mathematics and Statistics, University of Turku, 20014 Turku, Finland}
\email{mikko.m.jaskari@utu.fi}
\begin{document}
\begin{abstract}
We apply the resonance method to Montgomery's convolution formula for $\re\left(e^{-i\theta}\log\zeta(\sigma+it)\right)$ in the strip $1/2 < \sigma < 1$. This gives new insight into maximal values of $\re\left(e^{-i\theta}\log\zeta(\sigma+it)\right)$ for $t \in [T^{\beta},T]$ for all $\beta \in (0,1)$ and real $\theta$.
\end{abstract}
\maketitle
\section{Introduction}
The Riemann zeta function is a famously important function in number theory. All the non-trivial zeros of the zeta function are in the critical strip $0 \le \re(s) \le 1$ and they are related to the distribution of the prime numbers. This is one of the reasons why it is important to study the behavior of the zeta function inside the critical strip.

In 1977 H. L. Montgomery \cite{MR460255} proved that  if we let $\sigma \in (1/2,1)$ and $T > T_0(\sigma)$, then for any real $\theta$ we have\footnote{We denote $\log_{j}$ as the $j$th iterated logarithm and for instance $\log_{2}{T}=\log\log{T}$.}
\begin{align}
\max_{T^{(\sigma-1/2)/3} \le t \le T} \re\left( e^{-i\theta}\log\zeta(\sigma+it) \right) \ge \frac{1}{20}\left( \sigma-\frac{1}{2} \right)^{1/2}\frac{(\log{T})^{1-\sigma}}{(\log_{2}{T})^{\sigma}}. \label{eq1}
\end{align}
Moreover, assuming the Riemann hypothesis,  Montgomery \cite{MR460255} obtains
\begin{align}
\max_{T^{1/6} \le t \le T}\re\left( e^{-i\theta}\log\zeta(\sigma+it) \right) \ge \frac{1}{20}\frac{(\log{T})^{1-\sigma}}{(\log_{2}{T})^{\sigma}} \label{eq2}
\end{align}
for any real $\theta$. In Montgomery's method the lower bound for $t$ weakens notably without the use of the Riemann hypothesis when $\sigma \to 1/2^{+}$. C. Aistleitner improved the lower bound of the extreme value in \cite{MR3498919} in the case $\theta \equiv 0 \textrm{ (mod } 2\pi)$ by showing that for fixed $\sigma \in (1/2,1)$ and sufficiently large $T$ we have
\begin{align}
\max_{ 0 \le t \le T} |\zeta(\sigma+it)| \ge \exp\left( 0.18(2\sigma-1)^{1-\sigma}\frac{(\log{T})^{1-\sigma}}{(\log_{2}{T})^{\sigma}} \right).
\end{align} 
Our goal is to prove a Montgomery type result with better lower bound for $t$ when $\sigma$ is close to $1/2$.
\begin{theorem}\label{th1.1}
Fix $\sigma \in (1/2,1)$, $\beta \in (0,1)$ and ${0 < \kappa < \min(\sigma-1/2,1-\beta)}$. Then there exists a function $\upsilon : (1/2,1) \to \R_+$ which satisfies
\begin{align}
\upsilon(\sigma) &> c\min\left\{1,(e^2-e)\frac{|\log{(2\sigma-1)}|}{|\log{(2\sigma-1)}|+1}\right\},  \label{C2}
\end{align}
where $c$ is a positive constant independent from any of the chosen parameters $\sigma, \beta$ or $\kappa$ such that for any $\theta \in \R$ and sufficiently large $T$ we have 
\begin{align}
\max_{t \in [T^{\beta}, T]} \re\left( e^{-i\theta}\log\zeta(\sigma+it) \right) \ge \upsilon(\sigma)\frac{\kappa^{1-\sigma}}{\sqrt{|\log{(2\sigma-1)}|}}\frac{(\log{T})^{1-\sigma}}{(\log_{2}{T})^{\sigma}}.
\end{align}
\end{theorem}
Under the Riemann hypothesis we could choose $\kappa < 1 - \beta$ and would not require $\kappa < \sigma - 1/2$. Due to condition $\kappa < \sigma - 1/2$, our unconditional result in the case $t \in [0,T]$ is weaker than Montgomery's. Taking $\theta = \pi$ we obtain the following corollary.
\begin{corollary}\label{co1.1}
Fix $\sigma \in (1/2,1)$, $\beta \in (0,1)$ and ${0 < \kappa < \min(\sigma-1/2,1-\beta)}$. Then there exists a function $\upsilon(\sigma)$ satisfying the properties stated in Theorem \ref{th1.1} so that 
\begin{align}
\max_{t \in [T^{\beta}, T]} -\log|\zeta(\sigma+it)| \ge \upsilon(\sigma)\frac{\kappa^{1-\sigma}}{\sqrt{|\log{(2\sigma-1)}|}}\frac{(\log{T})^{1-\sigma}}{(\log_{2}{T})^{\sigma}}.
\end{align}
\end{corollary}
Result like corollary \ref{co1.1} can then be converted into estimate of the upper bound of the minimum of $|\zeta(\sigma+it)|$ for $t \in [T^{\beta}, T]$.

We will prove Theorem \ref{th1.1} by means of resonance method introduced by K. Soundararajan \cite{MR2425151}. Our work is heavily inspired by the works of A. Bondarenko and K. Seip \cite{MR3662441}, \cite{MR3880290} and \cite{MR3792093} and also by the work of A. Chirre and K. Mahatab  \cite{MR4379562}.

In \cite{MR3792093} Bondarenko and Seip proved that there exists a positive and continous function $\nu(\sigma)$ on $(1/2,1)$ bounded from below by $1/(2-2\sigma)$, with the asymptotic behavior
\begin{align}
\nu(\sigma) = \begin{cases} (1-\sigma)^{-1} + O(|\log{(1-\sigma)}|), & \sigma \to 1^{-} \\ (1/\sqrt{2} + o(1))\sqrt{|\log{(2\sigma-1)}|}, & \sigma \to 1/2^{+} \end{cases}
\end{align}
and such that the following holds. If $T$ is sufficiently large, then for $1/2 + 1/\log_{2}{T} \le \sigma \le 3/4$,
\begin{align}
\max_{t \in [\sqrt{T},T]} |\zeta(\sigma+it)| \ge \exp\left( \nu(\sigma)\frac{(\log{T})^{1-\sigma}}{(\log_{2}{T})^{\sigma}} \right)
\end{align}
and for $3/4 \le \sigma \le 1-1/\log_{2}{T}$,
\begin{align}
\max_{t \in [T/2,T]} |\zeta(\sigma+it)| \ge \log_{2}{T}\exp\left( c+\nu(\sigma)\frac{(\log{T})^{1-\sigma}}{(\log_{2}{T})^{\sigma}} \right)
\end{align}
with $c$ an absolute constant independent of $T$.

This theorem already improves both Aistleitner's and Montgomery's theorems and is seemingly much stronger than Theorem \ref{th1.1} in the specific case $\theta \equiv 0 $ (mod $2\pi$). The method Chirre and Mahatab used in \cite{MR4379562} to prove under Riemann hypothesis that for each fixed $\beta \in (0,1)$ there exists a positive constant $c > 0$ such that for sufficiently large $T$ we have
\begin{align}
\max_{ t\in [T^{\beta},T]} \pm\arg\zeta\left( \frac{1}{2} + it \right) \ge c\sqrt{\frac{\log{T}\log_{3}{T}}{\log_{2}{T}}}.
\end{align}
would also give stronger result to the case $\theta \equiv \pm\pi/2 $ (mod $2\pi$).  However, these proofs do not directly generalize to the case $\theta \in \R$ and particularly $\theta \equiv \pi$ (mod $2\pi$). 

Our strategy is to use Montgomery's convolution formula and apply resonance method by choosing a similar resonator as Bondarenko and Seip with some modifications. It is worth to mention that the result of Bondarenko and Seip has been recently improved by Zikang Dong and Bin Wei \cite{dong2022large} by a factor $2^{\sigma}$ as $\sigma \to 1/2^{+}$. Their result gives that if we let $\frac{1}{2} < \sigma < 1$ and fix $\beta \in (0,1)$ we have that
\begin{align}
\max_{t \in [T^{\beta}, T]} \log|\zeta(\sigma+it)| \ge c_{\beta}(\sigma)\frac{(\log{T})^{1-\sigma}}{(\log_{2}{T})^{\sigma}}
\end{align}
holds for a function $c_{\beta}(\sigma)$ which has the asymptotic behavior
\begin{align}
c_{\beta}(\sigma) = (\sqrt{2} +o(1))(1-\beta)^{1-\sigma}\sqrt{|\log{(\sigma-1/2)}|}, \textrm{ as } \sigma \to 1/2^{+}.
\end{align}
\section{The Convolution Formula}
We first introduce the needed convolution formula. The following lemma is \cite[Lemma 4]{MR460255}.
\begin{lemma}\label{lemma2.1}
Let $\sigma \in (1/2,1)$ and $t \ge 15$. Suppose that $\zeta(\sigma_0+iu) \not=0$ for any $\sigma_0 \in [\sigma,1)$ and any $u \in \R$ such that $|u-t| \le 2(\log{t})^2$. Then, for any $\psi > 0$ and any real $H$,
\begin{align}
\begin{split}
&\frac{2}{\pi}\int_{-(\log{t})^{2}}^{(\log{t})^{2}} \log\zeta(\sigma+i(t+u))\left( \frac{\sin{\psi u}}{u} \right)^2 e^{iHu}du \\ &= \sum_{n=1}^{\infty} \frac{\Lambda(n)\max(0,\psi-|H-\log{n}|)}{n^{\sigma+it}\log{n}} + O\left( \frac{e^{|H|+2\psi}}{(\log{t})^{2}} \right). \label{eq**}
\end{split}
\end{align}
\end{lemma}
Here $\Lambda$ is the von Mangoldt function. We now argue as Montgomery did in \cite[Equation (16)]{MR460255}. Let $\psi = 1/2$ and take successively $H = H_1,H_2,H_3$ where $H_1 := -\log{x}, H_2 := 0, H_3 := \log{x}$ where $1 \le x \le (\log{t})^2 $. Now the main term on the right hand side of (\ref{eq**}) vanishes when $H \in \{H_1, H_2\}$. Nevertheless, multiplying (\ref{eq**}) for $H_1, H_2$ and $H_3$ by $\frac{1}{2}e^{-i\theta}$, $1$ and $\frac{1}{2}e^{i\theta}$ respectively and adding them together we get
\begin{align}
\begin{split}
&\frac{2}{\pi}\int_{-(\log{t})^{2}}^{(\log{t})^{2}}\log\zeta(\sigma+i(t+u))\left( \frac{\sin{(u/2)}}{u} \right)^2 (1+\cos(\theta+u\log{x}))du \\
&= \frac{1}{2}  e^{i\theta}\sum_{e^{-1/2}x\le n \le e^{1/2}x} \frac{\Lambda(n)}{n^{\sigma+it}\log{n}}\left( \frac{1}{2} - \left| \log\frac{n}{x} \right| \right) + O\left( \frac{x}{(\log{t})^2} \right). \label{eq7}
\end{split}
\end{align}
\section{Proof of Theorem \ref{th1.1} by the Resonance Method}
The resonance method introduced by Soundararajan \cite{MR2425151} is based on evaluation of the ratio of two moments. The numerator moment is the integral of the investigated function multiplied by a chosen non-negative resonator over a chosen interval and the denominator moment is the integral of the resonator over the same interval. The ratio then gives a lower bound for the investigated function on that chosen interval.

In order to construct the resonator needed in the resonance method we will first define various sets motivated by \cite[Section 3]{MR3880290} and \cite[Section 4.4]{MR3792093}. Fix ${\sigma \in (1/2,1)}$, $\beta \in (0,1)$ and $0 < \kappa < \min(\sigma-1/2,1-\beta)$. Let us denote $N = \lfloor T^{\kappa} \rfloor$.
\begin{definition}\label{de3.1} Define the following sets \\
\begin{itemize}
\item $P$ is the set of primes in the interval $(e\log{N}\log_{2}{N},e^{2}\log{N}\log_{2}{N}].$ \\
\item Let $a > 1$ be fixed, let
\begin{align*} 
M := \left\{ m \in \N :  m \textrm{ has at least } \frac{a\log{N}}{|\log(2\sigma-1)|} \textrm{ prime factors in } P \right\}
\end{align*}
and let $M' \subset M$ contain all the integers in $M$ that have prime factors only in $P$.
\end{itemize}
\end{definition}
\begin{definition}\label{de3.2}
Let $f$ be the multiplicative function that is supported on the square-free numbers with all prime factors in $P$ and
\begin{align*}
f(p):= \frac{1}{\sqrt{|\log{(2\sigma-1)}|}}, \; p \in P. 
\end{align*}
\end{definition}
\begin{definition}\label{de3.3}
Define
\begin{align*}
\mathcal{M} := \textrm{supp}(f) \setminus  M
\end{align*}
which in other words means that
\begin{align*}
\mathcal{M} = \left\{ m \in \N : \begin{matrix} m \textrm{ is square-free, all prime factors of } m \textrm{ are in } P \\  \textrm{ and } m \textrm{ has at most } \frac{a\log{N}}{|\log(2\sigma-1)|} \textrm{ prime factors } \end{matrix} \right\}.
\end{align*}
Now let $\mathcal{J}$ be defined as
\begin{align*}
\mathcal{J} := \left\{ j \in \Z : [(1+T^{-1})^{j}, (1+T^{-1})^{j+1}) \cap \mathcal{M} \not= \emptyset \right\}
\end{align*}
and
\begin{align*}
\mathcal{M'} := \{ m_j :  j \in \mathcal{J}, m_j = \min\{[(1+T^{-1})^{j}, (1+T^{-1})^{j+1}) \cap \mathcal{M}  \} \}.
\end{align*}
\end{definition}
\begin{definition}\label{de3.4}
Let $\gamma \in (0,1)$ be fixed after $\sigma$ is fixed. Let $L$ be the set of integers that have at most $\frac{\gamma\log{N}}{|\log{(2\sigma-1)}|}$ prime factors in $P$. Let $L' \subset L$ contain integers in $L$ that have prime factors only in $P$. Now set
\begin{align*}
\mathcal{L} := \mathcal{M} \setminus L.
\end{align*}
Hence $\mathcal{L}$ is the set of integers in $\mathcal{M}$ that have at least $\frac{\gamma\log{N}}{|\log{(2\sigma-1)}|}$ prime factors.
\end{definition}
We now define the resonator we need.
\begin{definition}\label{de3.5}
Let $r : \mathcal{M}' \to \R$ be defined as
\begin{align*}
r(m_j) := \left( \sum_{\substack{(1-T^{-1})^{j-1} \le n \le (1+T^{-1})^{j+2} \\ n \in \mathcal{M}}} f(n)^2 \right)^{1/2}, \textrm{ for every } j \in \mathcal{J}
\end{align*}
and the resonator $R : \R \to \C$ as
\begin{align*}
R(t) := \sum_{m \in \mathcal{M}'} \frac{r(m)}{m^{it}}.
\end{align*}
\end{definition}
It is crucial that $|\mathcal{M}| \le N$ for large $N$ and this is essentially shown in \cite[Section 4.4]{MR3792093}. Note that the sets we use are subsets of the sets used in \cite[Section 4.4]{MR3792093} and we only consider the case $k=1$.

We will now move forward to the definition of the moments. First we denote
\begin{align}
\Phi(t) := e^{-t^2/2}
\end{align}
and let $\hat{\Phi}$ be the Fourier transform of $\Phi$ defined as
\begin{align}
\hat{\Phi}(y) := \int_{-\infty}^{\infty} \Phi(t)e^{-ity}dt = \sqrt{2\pi}\Phi(y). \label{eqFourier}
\end{align} 
Since in Lemma \ref{lemma2.1} we need to assume that there are no zeros of the Riemann zeta function on the right side of the contour we have to bypass all those cases where such zeros exists. We do so by defining the following indicator function.
\begin{definition}
Denote by $\rho$ non-trivial zeros of the $\zeta$-function and let $I(\sigma,t)$ be an indicator function defined as
\begin{align*}
I(\sigma,t) = \begin{cases} 1, &\textrm{if there is no zero } \rho \textrm{ such that } \re(\rho) \ge \sigma \textrm{ and } |t-\im(\rho)| \le (\log{t})^2 \\ 0, &\textrm{otherwise.} \end{cases}
\end{align*}
\end{definition}
Using the indicator function we define moments as follows.
\begin{definition}
Let $\sigma \in (1/2,1)$ be fixed. We define two moments $M_1(R,T)$ and $M_2(R,T)$ as 
\begin{align*}
M_1(R,T) &:= \int_{T^{\beta}}^{T\log{T}}\left( \int_{-(\log{t})^2}^{(\log{t})^2} K(u)du \right)|R(t)|^{2}\Phi\left( \frac{t}{T} \right)dt, \\
M_2(R,T) &:= \int_{T^{\beta}}^{T\log{T}}\left( \int_{-(\log{t})^2}^{(\log{t})^2} e^{-i\theta}\log\zeta(\sigma+i(t+u)) K(u)du \right)|R(t)|^{2}\Phi\left( \frac{t}{T} \right)I(\sigma,t)dt,
\end{align*}
where
\begin{align*}
K(u) &= \left( \frac{\sin{(u/2)}}{u} \right)^2 (1+\cos{(\theta+u\log{(e\log{N}\log_{2}{N})})}) \\ &+\left( \frac{\sin{(u/2)}}{u} \right)^2 (1+\cos{(\theta+u\log{(e^{3/2}\log{N}\log_{2}{N})})}) \\ &+\left( \frac{\sin{(u/2)}}{u} \right)^2 (1+\cos{(\theta+u\log{(e^{2}\log{N}\log_{2}{N})})}).
\end{align*}
\end{definition}
Now
\begin{align}
\max_{t \in [T^{\beta},T\log{T}]} \re(e^{-i\theta}\log\zeta(\sigma+it)) \ge \frac{\re(M_{2}(R,T))}{M_{1}(R,T)}. \label{eqResonance}
\end{align}
We define $P_{-}$ to be the set of primes in the interval $(e^{1/2}\log{N}\log_{2}{N},e\log{N}\log_{2}{N}]$ and $P_{+}$ to be the set of primes in the interval $(e^{2}\log{N}\log_{2}{N},e^{5/2}\log{N}\log_{2}{N}]$. We begin with the use of Lemma \ref{lemma2.1} and (\ref{eq7}) to obtain
\begin{align}
\begin{split}
& I(\sigma,t)\left( \int_{-(\log{t})^2}^{(\log{t})^2} e^{-i\theta}\log\zeta(\sigma+i(t+u)) K(u)du \right) \\ &= I(\sigma,t)\left(\frac{\pi}{4}\sum_{n \in P_{-}}\frac{\Lambda(n)}{n^{\sigma+it}\log{n}}\left( \frac{1}{2} - \left|\log{\frac{n}{e\log{N}\log_{2}{N}}} \right| \right)+\frac{\pi}{8}\sum_{n \in P}\frac{\Lambda(n)}{n^{\sigma+it}\log{n}} \right. \\ &+ \left.\frac{\pi}{4}\sum_{n \in P_{+}}\frac{\Lambda(n)}{n^{\sigma+it}\log{n}}\left( \frac{1}{2} - \left|\log{\frac{n}{e^{2}\log{N}\log_{2}{N}}} \right| \right) + O\left( 1  + \frac{\log_{2}{T}}{\log{T}} \right)\right)
\\ &= I(\sigma,t)\left(\frac{\pi}{4}\sum_{n \in P_{-}}\frac{\Lambda(n)}{n^{\sigma+it}\log{n}}\left( \frac{1}{2} - \left|\log{\frac{n}{e\log{N}\log_{2}{N}}} \right| \right)+\frac{\pi}{8}\sum_{n \in P}\frac{\Lambda(n)}{n^{\sigma+it}\log{n}} \right. \\ &+ \left.\frac{\pi}{4}\sum_{n \in P_{+}}\frac{\Lambda(n)}{n^{\sigma+it}\log{n}}\left( \frac{1}{2} - \left|\log{\frac{n}{e^{2}\log{N}\log_{2}{N}}} \right| \right) + O(1)\right)
\end{split} \label{eqN}
\end{align}
for $t \in [T^{\beta}, T\log{T}]$. In order to study the main term of $\re(M_{2}(R,T))$ we will first focus on the sum over $P$ in (\ref{eqN}) and deal with the sums over $P_{-}$ and $P_{+}$ later. Similarly to \cite[Section 4]{MR4379562} we obtain
\begin{align}
\begin{split}
&\re\left(\int_{T^{\beta}}^{T\log{T}}\frac{\pi}{8}\sum_{n \in P} \frac{\Lambda(n)}{n^{\sigma+it}\log{n}}|R(t)|^2\Phi\left( \frac{t}{T} \right)dt \right) \\ &= \frac{\pi}{8}\sum_{n \in P} \frac{\Lambda(n)}{n^{\sigma}\log{n}}\re\left( \int_{T^{\beta}}^{T\log{T}} n^{-it}|R(t)|^{2}\Phi\left( \frac{t}{T} \right) dt \right)\\
&= \frac{\pi}{8}\sum_{m,v \in \mathcal{M}'}\sum_{p \in P} \frac{r(m)r(v)}{p^{\sigma}} \re\left( \int_{T^{\beta}}^{T\log{T}} \Phi\left( \frac{t}{T} \right)e^{-it\log(mp/v)} dt\right).
\end{split} \label{eqM}
\end{align}
We first focus on evalution of
\begin{align}
\re\left( \int_{T^{\beta}}^{T\log{T}} \Phi\left( \frac{t}{T} \right)e^{-it\log(mp/v)} dt\right). \label{eq9}
\end{align}
Combining (\ref{eqFourier}) with the fact that $\Phi(t)$ is an even and real function we get
\begin{align}
\begin{split}
\re\left( \int_{0}^{\infty} \Phi\left( \frac{t}{T} \right)e^{-it\log(mp/v)} dt\right) &= \frac{1}{2}\int_{-\infty}^{\infty} \Phi\left( \frac{t}{T} \right)e^{-it\log(mp/v)}dt \\
&= T\frac{\sqrt{2\pi}}{2}\Phi\left( T\log\frac{mp}{v} \right). \label{eq10}
\end{split}
\end{align}
By the trivial estimate $|\Phi(u)e^{-iy}| \le 1$ we have
\begin{align}
\left| \re\left( \int_{0}^{T^{\beta}} \Phi\left( \frac{t}{T} \right)e^{-it\log(mp/v)} dt\right) \right| \le T^{\beta}. \label{eq11}
\end{align}
For $T > 193$ we have by rapid decay of $\Phi(t)$ as $t \to \infty$
\begin{align}
\left|\re\left( \int_{T\log{T}}^{\infty} \Phi\left( \frac{t}{T} \right)e^{-it\log(mp/v)} dt\right)\right| \le \int_{T\log{T}}^{\infty} \frac{1}{t^2} dt= o(1) \textrm{ as } T \to \infty. \label{eq12}
\end{align}
We may then conclude by (\ref{eq10}), (\ref{eq11}) and (\ref{eq12}) that since $\beta < 1$ we have
\begin{align}
\re\left( \int_{T^{\beta}}^{T\log{T}} \Phi\left( \frac{t}{T} \right)e^{-it\log(mp/v)} dt\right) = T\frac{\sqrt{2\pi}}{2}\Phi\left( T\log\frac{mp}{v} \right) + O(T^{\beta}).  \label{eqc}
\end{align}
In order to evaluate $M_{2}(R,T)$ we have to remove all $t \in [T^{\beta},T\log{T}]$ with ${I(\sigma,t)=0}$.
\begin{lemma}\label{lemma3.7}
Let $N(\sigma,T)$ denote the number of zeros $\rho$ of the $\zeta$-function for which $\re(\rho) \ge \sigma$ and $ 0 \le \im(\rho) \le T$. Now, for $T \ge 10$ and $1/2 \le \sigma \le 1$,
\begin{align*}
N(\sigma,T) \ll T^{3/2-\sigma}(\log{T})^{5}.
\end{align*}
\begin{proof} See \cite[Lemma 5]{MR460255}. \end{proof}
\end{lemma}
We note that $\textrm{meas}\{ t \in [T^{\beta},T\log{T}]: I(\sigma,t)=0\} \ll N(\sigma,T\log{T})(\log{T})^{2}$. By Lemma \ref{lemma3.7} and similar estimate as done for (\ref{eq11}) we may conclude that
\begin{align*}
\re\left( \int_{T^{\beta}}^{T\log{T}} \Phi\left( \frac{t}{T} \right)e^{-it\log(mp/v)} (1-I(\sigma,t))dt\right) \ll T^{3/2-\sigma}(\log{T})^{9}.
\end{align*}
Hence
\begin{align}
\begin{split}
&\re\left( \int_{T^{\beta}}^{T\log{T}} \Phi\left( \frac{t}{T} \right)e^{-it\log(mp/v)} I(\sigma,t)dt\right) \\ &= T\frac{\sqrt{2\pi}}{2}\Phi\left( T\log\frac{mp}{v} \right) + O(T^{3/2-\sigma}(\log{T})^{9} + T^{\beta}).
\end{split}  \label{eqA}
\end{align}
In order to evaluate the contribution of the error terms we have to estimate the size of the resonator.
\begin{lemma}\label{lemmaMid}
For any real $t$ we have
\begin{align*}
|R(t)|^2 \le 3T^{\kappa}\sum_{l\in \mathcal{M}}f(l)^2.
\end{align*}
\begin{proof}
We follow \cite[Proof of Lemma 5]{MR3880290}. Using $|R(t)| \le R(0)$ we begin with
\begin{align}
 R(0)^2 &= \sum_{m,n \in \mathcal{M}'} r(m)r(n) \le |\mathcal{M}'|\sum_{m \in \mathcal{M}'} r(m)^2 \label{eqXX}
\end{align}
where we used the inequality $ab \le (a^2 +b^2)/2$. Recall from the beginning of this section and notes after Definition \ref{de3.5} that there are at most $N = \lfloor T^{\kappa} \rfloor$ elements in $\mathcal{M}'$. Now by the definition of the set $\mathcal{M}'$ and the function $r$
\begin{align}
\sum_{m \in \mathcal{M}'} r(m)^2 &\le 3\sum_{l \in \mathcal{M}} f(l)^2. \label{eqYY}
\end{align}
We obtain the desired result by combining the upper bound of $|\mathcal{M}'|$, (\ref{eqXX}) and (\ref{eqYY}).
\end{proof}
\end{lemma}
We note that the sums over $P_{-}$ and $P_{+}$ in (\ref{eqN}) can be handled in the same way as sum over $P$ and thus we can ignore the non-negative contribution they give. Combining the definition of $M_{2}(R,T)$, (\ref{eqN}), (\ref{eqM}), (\ref{eqc}), (\ref{eqA}) and Lemma \ref{lemmaMid} we obtain
\begin{align}
\begin{split}
\re(M_{2}(R,T)) &\ge T\frac{\pi\sqrt{2\pi}}{16}\sum_{m,v \in \mathcal{M}'}\sum_{p \in P}\frac{r(m)r(v)}{p^{\sigma}}\left(\Phi\left( T\log\frac{mp}{v} \right) \right) \\ &+ O\left((T^{3/2+\kappa-\sigma}(\log{T})^{9} + T^{\beta+\kappa})\sum_{l\in\mathcal{M}}f(l)^2 \frac{(\log{N})^{1-\sigma}}{(\log_{2}{N})^{\sigma}}\right). \label{eqD}
\end{split}
\end{align}
The next two lemmas allow us to lower bound the main term on the right hand side of (\ref{eqD}).
\begin{lemma}\label{lemma3.9}
We have
\begin{align*}
\sum_{m,v \in \mathcal{M}'}\sum_{p \in P} \frac{r(m)r(v)}{p^{\sigma}}\Phi\left( T\log\frac{mp}{v} \right) \ge \sum_{v \in \mathcal{M}} f(v)^2 \sum_{p|v}\frac{1}{f(p)p^{\sigma}}.
\end{align*}
\begin{proof}
We follow \cite[Section 3]{MR3880290}. We consider all triples $m', v' \in \mathcal{M}'$ and $p \in P$ such that $|pm'/v'-1| \le3/T$. We use the notation 
\begin{align*}
J(m') := [(1+T^{-1})^{j}, (1+T^{-1})^{j+1}),
\end{align*}
where $j$ is the unique integer such that $(1+T^{-1})^{j} \le m < (1+T^{-1})^{j+1}$. By the definition of $r(m')$ and the Cauchy-Schwarz inequality we have for any $p \in P$ and any $m',v' \in \mathcal{M}'$
\begin{align*}
\sum_{\substack{m,v \in \mathcal{M} \\ mp=v \\ m \in J(m'), v \in J(v') }} f(m)f(v) &\le \left(\sum_{\substack{m,v \in \mathcal{M} \\ mp=v \\ m \in J(m'), v \in J(v') }} f(m)^2\right)^{1/2}\left(\sum_{\substack{m,v \in \mathcal{M} \\ mp=v \\ m \in J(m'), v \in J(v') }} f(v)^2\right)^{1/2} \\ &\le \left( \sum_{m \in J(m')} f(m)^2 \right)^{1/2}\left( \sum_{v \in J(v')} f(v)^2 \right)^{1/2} \\ &\le r(m')r(v')
\end{align*}
and hence, by the definition of $\mathcal{M}'$, that, for any $p \in P$,
\begin{align*}
\sum_{\substack{m,v \in \mathcal{M} \\ mp=v}} f(m)f(v) \le \sum_{\substack{m', v' \in \mathcal{M}' \\ |pm'/v'-1| \le3/T }} r(m')r(v').
\end{align*}
Now
\begin{align*}
\sum_{m,v \in \mathcal{M}'}\sum_{p \in P} \frac{r(m)r(v)}{p^{\sigma}}\Phi\left( T\log\frac{mp}{v} \right) \ge \sum_{p \in P}\sum_{\substack{m,v \in \mathcal{M} \\ mp=v}} \frac{f(m)f(v)}{p^{\sigma}} = \sum_{v \in \mathcal{M}} f(v)^2 \sum_{p|v}\frac{1}{f(p)p^{\sigma}}.
\end{align*}
In the last step we used multiplicativity of $f$. This completes the proof.
\end{proof}
\end{lemma}
\begin{lemma}\label{lemma3.12}
For $\gamma \in (0,1)$ we have
\begin{align*}
\sum_{v \in \mathcal{M}} f(v)^2 \sum_{p|v}\frac{1}{f(p)p^{\sigma}} \ge \gamma\sum_{v \in \mathcal{L}}f(v)^2\frac{ e^{-2\sigma}}{\sqrt{|\log{(2\sigma-1)}|}}\frac{(\log{N})^{1-\sigma}}{(\log_{2}{N})^{\sigma}}
\end{align*}
\begin{proof}
We follow similar ideas as in \cite[Proof of Lemma 4]{MR3880290}. We recall the definition of the set $\mathcal{L}$ from Definition \ref{de3.4}.
We have
\begin{align}
\sum_{v \in \mathcal{M}}f(v)^2\sum_{p|v} \frac{1}{f(p)p^{\sigma}} \ge \sum_{v \in \mathcal{L}}f(v)^2\sum_{p|v} \frac{1}{f(p)p^{\sigma}}. \label{eq16}
\end{align}
Now
\begin{align}
\sum_{v \in \mathcal{L}}f(v)^2\sum_{p|v} \frac{1}{f(p)p^{\sigma}} \ge \sum_{v \in \mathcal{L}}f(v)^2 \frac{\gamma\log{N}}{|\log{(2\sigma-1)}|}\min_{p \in P} \frac{1}{f(p)p^{\sigma}}. \label{eq17}
\end{align}
We obtain the desired result by noting that
\begin{align*}
\min_{p \in P} \frac{1}{f(p)p^{\sigma}} \ge \sqrt{|\log{(2\sigma-1)}|}(e^{2}\log{N}\log_{2}{N})^{-\sigma}
\end{align*}
and combining  (\ref{eq16}) and (\ref{eq17}).
\end{proof}
\end{lemma}
Now by combining (\ref{eqD}) and Lemmas \ref{lemma3.9} and \ref{lemma3.12} we obtain the following lemma.
\begin{lemma}\label{lemmaProp}
Fix $\sigma \in (1/2,1)$, $\beta \in (0,1)$ and $\kappa < \min(\sigma-1/2, 1-\beta)$. Then there exists a positive constant $c_1$ independent from $\sigma$, $\beta$ and $\kappa$ such that
\begin{align}
\begin{split}
\re(M_{2}(R,T)) \ge & c_{1}\gamma T\sum_{v \in \mathcal{L}} f(v)^{2}\frac{1}{\sqrt{|\log{(2\sigma-1)}|}}\frac{(\log{N})^{1-\sigma}}{(\log_{2}{N})^{\sigma}} \\ &+O\left((T^{3/2+\kappa-\sigma}(\log{T})^{9} + T^{\beta+\kappa})\sum_{v \in \mathcal{M}}f(v)^{2}\right).
\end{split} \label{eqMain2}
\end{align}
\end{lemma}
Here the restriction $\kappa < \min(\sigma-1/2,1-\beta)$ guarantees that the error term is acceptable. We now move forward to the evaluation of the denominator moment $M_1(R,T)$.
\begin{lemma}\label{lemmaM1}
There exists a positive constant $c_2$ such that
\begin{align*}
M_{1}(R,T)  \le c_2 T\sum_{n \in \mathcal{M}}f(n)^2.
\end{align*}
\begin{proof}
Since $K(u) \ge 0$ for all $u$ note that as in \cite[Proof of Theorem 1]{MR460255} we have
\begin{align}
\int_{-\infty}^{\infty}K(u)du \ll 1 \label{eql1}
\end{align}
since $K(u)$ is a sum of three distinct kernels introduced in Montgomery's paper. The proof then follows from the argument made in \cite[Proof of Lemma 5]{MR3880290}.
\end{proof}
\end{lemma}
Now according to (\ref{eqResonance}) and Lemmas \ref{lemmaProp} and \ref{lemmaM1} there is a need to evaluate the ratio
\begin{align*}
\frac{\sum_{v \in \mathcal{L}}f(v)^2}{\sum_{n\in\mathcal{M}}f(n)^2}
\end{align*}
and this is done in our next lemmas.
\begin{lemma}\label{lemmaMM}
Let $a>1$ be fixed and $M$ is defined as in Definition \ref{de3.1}. Then
\begin{align*}
\frac{1}{\sum_{j \in \N} f(j)^2} \sum_{v \in M} f(v)^2 = o(1) \quad \textrm{as } N \to \infty.
\end{align*}
\begin{proof}See \cite[Proof of Lemma 8]{MR3792093}.\end{proof}
\end{lemma}
\begin{lemma}\label{lemma5.8}
Let $\sigma \in (1/2,1)$, $b \in (0,1)$ and $\gamma$ which appears in Definition \ref{de3.4} for the set $\mathcal{L}$ satisfy 
\begin{align}
\gamma < \min\{1,C(\sigma)(e^2-e)(b-1)/\log{b}\} \label{gamma1}
\end{align} 
where
\begin{align*}
C(\sigma) := \frac{1}{1+f(p)^2}.
\end{align*}
Then
\begin{align*}
\frac{1}{\sum_{j \in \N} f(j)^2}\sum_{v \not\in \mathcal{L}} f(v)^2 = o(1) \quad \textrm{as } N \to \infty.
\end{align*}
\begin{proof}
We follow the ideas in \cite[Proof of Lemma 4]{MR3880290}. We note that
\begin{align*}
\mathcal{L} = \textrm{supp}(f) \setminus (M \cup L)
\end{align*}
and by Lemma \ref{lemmaMM} it is sufficient to show
\begin{align*}
\frac{1}{\sum_{j \in \N} f(j)^2}\sum_{v \in L} f(v)^2 = o(1) \quad \textrm{as } N \to \infty.
\end{align*}
We recall the definition of $L'$ from Definition \ref{de3.4} and use the fact that $f$ is multiplicative to obtain
\begin{align*}
\frac{1}{\sum_{j \in \N} f(j)^2}\sum_{v \in L} f(v)^2 = \frac{1}{\prod_{p \in P} (1+f(p)^2)}\sum_{v \in L'}f(v)^2
\end{align*}
where, for any $b \in(0,1)$,
\begin{align*}
\sum_{v \in L'} f(v)^2 \le b^{-\gamma\frac{\log{N}}{|\log{(2\sigma-1)}|}}\prod_{p \in P}(1+bf(p)^2).
\end{align*}
For $b$ that is chosen close enough to $1$, we have
\begin{align*}
\frac{1}{\prod_{p \in P} (1+f(p)^2)}\sum_{v \in L'}f(v)^2 \le b^{-\gamma\frac{\log{N}}{|\log{(2\sigma-1)}|}}\exp\left( C(\sigma)\sum_{p \in P} (b-1)f(p)^2 \right),
\end{align*}
due to the Taylor expansion of 
\begin{align*}
\log{\left( \frac{1+bf(p)^2}{1+f(p)^2} \right)}=\log{\left(\frac{1-(1-b)f(p)^2}{1+f(p)^2}\right)}. 
\end{align*}
The cardinality of $P$ is at most $e^2\log{N}$.  By the prime number theorem
\begin{align*}
\sum_{p \in P} f(p)^2 = (1+o(1))\frac{\log{N}}{|\log{(2\sigma-1)}|}(e^2 -e).
\end{align*}
Hence we obtain
\begin{align*}
&\frac{1}{\sum_{j \in \N} f(j)^2}\sum_{v \in L} f(v)^2 \\ &\le \exp\left(\left( C(\sigma)(e^2-e)(b-1)-\gamma\log{b}+o(1) \right)\frac{\log{N}}{|\log{(2\sigma-1)}|} \right).
\end{align*}
If we choose $b$ close to $1$ we may note that
\begin{align*}
C(\sigma)(e^2-e)(b-1) - \gamma\log{b} < 0
\end{align*}
holds for $\gamma$ satisfying (\ref{gamma1}). This completes the proof.
\end{proof}
\end{lemma}
Since $\mathcal{L} = \textrm{supp}(f) \setminus (M \cup L)$ and $\mathcal{M} = \textrm{supp}(f) \setminus M$ we have by Lemmas \ref{lemmaMM} and \ref{lemma5.8}
\begin{align}
\frac{\sum_{v \in \mathcal{L}}f(v)^2}{\sum_{n\in\mathcal{M}}f(n)^2} = 1-o(1), \quad \textrm{as } N \to \infty. \label{eq31}
\end{align}
We note that in Lemma \ref{lemma5.8} if we choose largest possible $\gamma$, the lower bound in Lemma \ref{lemmaProp} is also largest possible and we can choose $\gamma$ close to $1$ if $\sigma \le 0.880766...$ holds. Combining this observation with (\ref{eqResonance}), (\ref{eq31}), Lemmas \ref{lemmaProp} and \ref{lemmaM1} and restriction (\ref{gamma1}) we may conclude that
\begin{align}
\begin{split}
&\max_{t \in [T^{\beta},T\log{T}]} \re(e^{-i\theta}\log\zeta(\sigma+it)) \\ &\ge  \upsilon(\sigma)\frac{\kappa^{1-\sigma}}{\sqrt{|\log{(2\sigma-1)}|}}\frac{(\log{T})^{1-\sigma}}{(\log_{2}{T})^{\sigma}}+O\left(T^{1/2+\kappa-\sigma}(\log{T})^{9} + T^{\beta+\kappa-1}\right). \label{eq32}
\end{split}
\end{align}
for $\upsilon : (1/2,1) \to \R_+$ satisfying (\ref{C2}). Theorem \ref{th1.1} then follows by noting that the desired restriction $T^{\beta} \le t \le T$ is obtained by trivial adjustment, applying (\ref{eq32}) for $T/\log{T}$ in place of $T$ and $\beta' \in (\beta,1-\kappa)$ in place of $\beta$.
\section*{Acknowledgements}
I was partially funded by UTUGS graduate school and completed this work while working in Academy of Finland projects no. 333707 and 346307. I would also like to thank Kaisa Matomäki for supervising my work, Kamalakshya Mahatab for giving ideas for this project and the anonymous referee for valuable comments.
\bibliographystyle{plain}
\bibliography{mybib}
\end{document}